\documentclass{article}

\usepackage{amsmath,amsthm,amsfonts}
\usepackage{amssymb}
\usepackage{graphicx}
\usepackage{subfigure}
\usepackage{rotate}

\def\sl R{\mathscr{R}}

\def\eps{\epsilon}

\begin{document}


\begin{center}
{\Large\bf Monte Carlo without Chains}
\end{center}
\vskip14pt

\begin{center}
{\bf Alexandre J.\ Chorin}
\end{center}
\vskip12pt
\centerline{Department of Mathematics, University of California}
\centerline{and}
\centerline{Lawrence Berkeley National Laboratory}
\centerline{Berkeley, CA 94720}

\begin{abstract}
A sampling method for spin systems is presented.
The spin lattice is written as the union of a nested sequence of sublattices, all but the last with conditionally
independent spins, which are sampled in succession using their marginals.
The marginals are computed concurrently by a fast
algorithm;
errors in the evaluation of the marginals are offset by weights.
There are no Markov chains and each sample is independent of the previous ones;
the cost of a sample is proportional to the number of spins
(but the number of samples needed for good statistics may grow with array size).
The examples include the Edwards-Anderson spin glass in three dimensions.
\end{abstract}

\section{Introduction.}

Monte Carlo sampling in physics is synonymous with Markov chain Monte Carlo
(MCMC) for good reasons which are too well-known to need repeating (see e.g. \cite{Bi2},\cite{Li1}).
Yet there are problems where the free energy landscape exhibits multiple
minima and the number of MCMC steps needed to produce an independent sample
is huge; this happens for example in spin glass models (see e.g. \cite{Fi1},\cite{Ka3}).
It is therefore worthwhile to consider alternatives, and the
purpose of the present paper is to propose one.

An overview of the proposal is as follows:
Consider a set of variables (``spins") located at the nodes of a lattice $L$
with a probability density $P$ that one wishes to sample.
Suppose one can construct a nested sequences of subsets
$L_0 \supset L_1 \supset \dots \supset L_n$
with the following properties: $L_0=L$; $L_n$ contains few points; the marginal density of the variables in each $L_i$ is known, and given values of the spins in $L_{i+1}$, the remaining variables in $L_i$ are independent.
Then the following is an effective sampling strategy for
the spins in $L$: First sample the spins in $L_n$ so that each configuration is sampled with a
frequency equal to its probability by first listing all the states of the spins in $L_n$ and calculating their probabilities. Then sample the variables in each $L_{i-1}$
as $i$ decreases from $n-1$ to zero using the independence of these 
variables, making sure that each state is visited with a frequency equal to its probability. Each state of $L=L_0$ is then also sampled
with a frequency equal to its probability, achieving importance sampling; the cost of each
sample of $L$ is proportional to the number of spins in $L$, and two successive samples are independent.
This can be done exactly in a few uninteresting cases (for example in the one-dimensional Ising model, see e.g. \cite{Ka1}),
but we will show that it can often be done approximately. The errors that come from the approximation can then be compensated for through the use of sampling weights.
The examples shown below include
spin glass models.

The heart of the construction is the fast evaluation of marginals, of which a previous version was presented in \cite{Ch9}; this is the subject of section 2. An example of the nested sets needed in the construction above is presented in section 3.
The construction is in general approximate, and the resulting 
errors are to be compensated for by weights, whose calculation is explained in section
4; efficiency demands a balance between the accuracy of the marginalization and
the variability of the weights, as explained in sections 4 and 6. The fast evaluation of marginals requires a sampling algorithm but the sampling
algorithm
is defined only once the marginalization is in place; this conundrum is resolved by
an iteration which is presented in section 5; two iteration steps turn out
to be sufficient. 

The algorithm is applied to the two-dimensional Ising model in section 6;
this is a check of self-consistency. Aspects of the Edwards-Anderson  (EA) spin
glass model three dimensions are discussed in section 7. Extensions and conclusions are presented in a concluding section.

As far as I know, there is no previous published work on Monte Carlo without chains for
spin systems. The construction of marginals explained below and in the earlier paper \cite{Ch9} is
a renormalization in the  sense of Kadanoff (more precisely, a decimation) \cite{Ka1};
renormalization has been used by many authors as a tool for speeding up
MCMC, see e.g. \cite{Bi3},\cite{Br1},\cite{Go4} . A construction conceptually related to the one here and also based on marginalization, but with a Markov chain, was presented in \cite{We1}.
An alternative construction of marginals can be found in \cite{St1}. There is some kinship between the construction here and the
decimation and message passing constructions in \cite{Bra1},\cite{Co1}. The specific connection between marginalization and conditional expectation used here originated in work on system reduction in the
framework of optimal prediction, see \cite{Chth},\cite{Ch13}.

The results in this paper are preliminary in the sense that the marginalization is performed in the simplest way I could imagine; more sophisticated versions
are suggested in the results and conclusion sections. Real-space
renormalization or decimation and the evaluation of marginals are one and the same, and the present
work could be written in either physics or probability language; it is the second
option that has been adopted.

All the examples below are of spin systems with two-valued spins and near-neighbor interactions in either a square lattice in two dimensions or a cubic lattice in three dimensions. The
Hamiltonians have the form $H=-\sum s_{i,j,k}\sum' J_{i,j,k,\ell}s_{i',j',k'}$ (with one subscript less in two dimensions), where the summation $\sum'$ is over near neighbors. In the Ising case the $J_{i,j,k,\ell}$ are independent of the indices, in the spin glass case $J_{i,j,k,\ell}$ are independent random variables; the index $\ell$ labels the direction of the interaction.


\section{Fast evaluation of approximate marginals.}

In this section we present an algorithm for the evaluation of marginals, which
is an updated version of the algorithm presented in \cite{Ch9}. For simplicity, we assume in this section, as in \cite{Ch9}, that the spins live on a two-dimensional
lattice with periodic boundary conditions (the generalization to three
dimensions is straightforward
except for the geometrical issues discussed in section 3).
We show how to go from $L_0$, the set of spins at the points $(i,j)$ of a regular
lattice whose probability density is known, to $L_1$, set of spins at the
points such that $(i+j)$ is even, and whose marginal is sought.
The probability density function (pdf) of the variables in $L_1$ can be written in the form $P_0=e^{W^{(0)}}/Z$, where
$W^{(0)}=-\beta H$, $\beta$ is the inverse temperature, $Z$ is the normalization constant, and the Hamiltonian $H$ has the form specificed in the introduction. To simplify notations, write $\tilde J_{i,j,\ell}=-\beta J_{i,j,\ell}$ and then drop the tildes,  so that
\begin{equation}
W^{(0)}=\sum
s_{i,j}\left(J_{i,j,1}s_{i+1,j}+J_{i,j,2} s_{i,j+1}\right).
\label{Hamiltonian}
\end{equation}
The $s_{i,j}$ take the values $\pm 1$.

Let the marginal density of the spins in $L_1$ be $P_1$. One can always write
$P_1=e^{W^{(1)}}/Z$, where $Z$ is the same constant as in the pdf of $P_0$.
Call the set of spins
in $L_1$ ``$\hat S$",
and the set of spins in $L_0$ but not in $L_1$ ``$\tilde S$",
so that $S=\hat S \cup \tilde S$  is the set of spins in $L_0$.
By definition of a marginal,
\begin{equation*}
P_1(\hat S)=e^{W^{(1)}}=\sum_{\tilde S} e^{W^{(0)}(S)}
\end{equation*}
or
\begin{equation}
W^{(1)}=\log \sum_{\tilde S} e^{W^{(0)}(S)},
\label{okunev}
\end{equation}
where the summation is over all the values of the spins in $\tilde S$.
Extend the range of values of the spins in $\hat S$ (but not
$\tilde S$) to the interval [0,1] as continuous variables but leave the expression for
the Hamiltonian unchanged (this device
is due to Okunev \cite{0k1} and replaces the more awkward construction in \cite{Ch9}).
Differentiate equation (\ref{okunev})  with respect to one of the a newly continuous variables $s=s_{i,j}$
(we omit the indices $(i,j)$ to make the formulas easier to read):
\begin{equation*}
\frac{\partial W^{(1)}}{\partial s}=W^{(1)'}=
\frac{
\sum_{\tilde S}
\frac{\partial W^{(0)}}{\partial s} e^{W^{(0)}(S)}}
{\sum_{\tilde S} e^{W^{(0)}(S)}},
\end{equation*}
or
\begin{equation}
\frac{\partial W^{(1)}}{\partial s}=E\left[\frac{\partial W^{(0)}}{\partial s}\mid\hat S\right],
\label{cond}
\end{equation}
where $E[\cdot\mid\hat S]$ denotes a conditional expectation given
$\hat S$.
A conditional expectation given $\hat S$ is an orthogonal projection onto the
space of functions of $\hat S$, and we approximate it by projecting onto the
span of a finite basis of functions of $\hat S$.

Before carrying out this projection, one should note the following property of $W^{(1)'}=\frac{\partial W^{(1)}}{\partial s}$:
take two groups of spins distant from each other in space, say $S_1$ and $S_2$.
The variables in these
groups should be approximately independent of each other, so that their joint pdf is
approximately the product of their separate pdfs.
The logarithm of their joint pdf is approximately the sum of the logarithms of their separate pdfs, and the derivative of that logarithm with respect to a variable
in $S_1$ should not be a function of the variables in $S_2$.
As a result, if one expands $\frac{\partial W}{\partial s}$ at
$s=s_{i,j}$, one needs only to project on a set of functions of $s_{i,j}$
and of a few neighbors of the point $(i,j)$. It is this observation that makes the algorithm in the present section
effective, and it is implied in the Kadanoff construction of a renormalized Hamiltonian, see e.g. \cite{Ka1}.

As a basis on which to project, consider, following Kadanoff, the
polynomials in $\hat{S}$ of the form:
$\psi_{p,q}=\sum_{i,j}s_{i,j}s_{i+p,j+q}$ for various values
of $p,q$, as well as polynomials of higher degree in the variables
$\hat{S}$. Define $\psi'_{p,q}=\partial \psi_{p,q}/\partial s_{i,j}$; the functions $\psi'_{p,q}$ involve only near neighbors of $s_{i,j}$ (for example, if $\psi_{1,1}=\sum_{\ell,k} s_{\ell,k}s_{\ell+1,k+1}$,
then $\psi'_{1,1}=2(s_{i+1,j+1}+s_{i-1,j-1})$ (no summation). Write the approximate conditional expectation of $W^{(1)'}$  as a sum: \begin{equation}
E[W^{(1)'}\mid\hat S]=\sum a_{p,q}\psi_{p,q}^{'}+ \cdots.
\label{expansion}
\end{equation}
Each function $\psi_{p,q}$ embodies an interaction, or linkage, between spins $(p,q)$ apart, and this is an expansion in ``successive linkages". The functions $W^{(0)}, W^{(1)},$ are invariant under the global symmetry $s \rightarrow  -s$,
and only polynomials having this symmetry need to be considered
(but see the discussion of symmetry breaking in section 6). For the reasons stated above, this series should converge rapidly as $p,q$ increase.
Evaluate the coefficients in (\ref{expansion})
by orthogonal projection onto the span of the $\psi'_{p,q}$.
This produces one equation per point $(i,j)$ (unless the system
is translation invariant, like the Ising model, in which case all these
equations are translates of each other). Assume furthermore that one has an algorithm for sampling the pdf
$P_0=e^{W^{(0)}}/Z$ (this is not trivial as the goal of the whole
exercise is to find good ways to sample $P_0$;
see section 5).
The projection can then be carried out by the usual method:
reindex the basis functions with a single integer, so that they become
$\psi_1, \psi_2, \ldots$, say $\psi_1=\psi_{1,1}$ etc.; at each point $(i,j)$ estimate by Monte Carlo the entries $a_{p,q}=E[\psi_p \psi_q]$ of a matrix $A$,
and the entries $b_p=E[W^{(0)'}\psi_p]$ of a vector $b$, where $E[\cdot]$ denotes
an expectation.
The projection we want is $\sum a_p\psi'_p$, where the
coefficients $a_p$
are the entries of the vector $A^{-1}b$ (see e.g. \cite{Ch13}).
In the current paper we use only the simplest basis with
$\psi_{1,1}, \psi_{1,-1}, \psi_{-1,1}, \psi_{-1,-1}$ (functions such as $\psi_{0,1}$ or $\psi_{1,0}$ do not appear because they involve spins not in $\hat S$).
In three
dimensions also we use basis functions of the form
$\sum s_{i,jk} s_{i',j'k'}$ where $(i',j',k')$ is a near neighbor of $(i,j,k)$
on a reduced lattice.
The locality of the functions $\psi'$ make the algorithm efficient;
more elaborate bases for the Ising case can be found in \cite{Ch9}. We have not invoked here any translation invariance, in view of later
applications to spin glasses.
For the Hamiltonian (\ref{Hamiltonian}), the quantity $W^{(0)'}=\frac{\partial W^{(0)}}{\partial s}$ for $s=s_{i,j}$ is
\begin{equation} \frac{\partial W^{(0)}}{\partial s}\mid_{i,j}=J_{i,j,1}s_{i+1,j} + J_{i-1,j,1}s_{i-1,j} +J_{i,j,2}s_{i,j+1} + J_{i,j-1,2} s_{i,j-1}.
\end{equation} The marginal density of the variables in $L_2$ (the set of points $(i,j)$ such that both $i$ and $j$ are odd),
is obtained by projecting on the basis functions
$\psi_{0,2}, \psi_{0,-2},$ $\psi_{2,0}, \psi_{-2,0}$, etc.
A single sample of all the spins in $L_0$ can be used  to generate samples of the inner products
needed to find the coefficients in the projections on all the sublattices; it is not a good idea to use the marginal of $L_1$ to evaluate the marginal of $L_2$, etc.,
because this may lead to a catastrophic error accumulation \cite{Ch9},\cite{Sw2};
it is the original $W{^{0'}}=\frac{\partial W^{(0)}}{\partial s}$ that is projected in the expansions at the different levels.

The last step is to reconstruct $W^{(1)}$ from its derivatives $W^{(1)'}$. In the
Ising (translation invariant) case, this is trivial: $W^{(1)'}=\sum a_p\psi_p^{'}$
implies $W^{(1)}=\sum a_p \psi_p$. In the spin glass case, there is
a minor conceptual (thought not practical) difficulty. For the various computed functions $W^{(1)'}$ to
be derivatives of some $W^{(1)}$ their cross derivatives must be equal. However, this equality requires an equality between coefficients evaluated at different points $(i,j)$ of the lattice. For example,
if $\psi_1$ after renumbering is $\psi_{1,1}$ before the renumbering in the notations above, and $\psi_3$ is $\psi_{-1,-1}$ before the renumbering, then the coefficient $a_3$ at the point $(i,j)$ should equal
the coefficient $a_1$ at the point $(i-1,j-1)$, and indeed these coefficients
describe the same interaction between the spins at $(i,j)$ and $(i-1,j-1)$.
However, these coefficients are computed separately by approximate computations,
and therefore, though closely correlated (with a correlation coefficient
typically above $.95$), they are not identical. This is not a practical
difficulty, because the replacement of one of these coefficients by the other,
or of both by some convex linear combination of the two, does not
measurably affect the outcome of the calculation. The conceptual issue is resolved if one notices that (i) for every lattice
$L_j$, except the smallest one $L_n$, one needs the coefficients $a_p$ at only half
the points, and using the coefficients calculated at these points is unambiguous, and (ii)
allowing these coefficients to differ is the same as writing the Hamiltonian
$W$ as $(s,Ms)$ where $s$ is the vector of spins and $M$ is an asymmetric
matrix. However, the values of $W$ are the same if $M$ is replaced by its symmetric part, which means replacing each of the two values of a
coefficient by the mean of the two values. If this is done everywhere
the cross derivatives become equal.

Finally, a practical comment. One may worry about a possible loss of
accuracy due to the ill-conditioning of a projection on a non-orthogonal polynomial
basis. I did develop an approximate Monte Carlo
Gram-Schmidt orthogonalization
algorithm and compared the resulting projection with
what has just been described; I could see no difference.


\section
{The nested sublattices.}

In this section we construct nested sequences of sublattices
$L=L_0, L_1,\dots$ such that the spins in $L_j$ are independent once
those in $L_{j+1}$ are determined. There is nothing unique about this construction; the nested sequences
should have appropriate independence and approximation properties while leading to
efficient programs.

The two-dimensional case was already discussed in the previous section: assume
$L=L_0$ is the set of nodes $(i,j), i,j,$ integers. We can choose as the next smaller array of spins $L_1$ the set of spins
at the location $(i,j)$ with $i+j$ even; then the spins in $L_0$ are independent
once in $L_1$ are known. The next lattice $L_2$ is the one where $i,j$ are both odd;
if the marginal on $L_1$ is approximated  with the four basis functions described in the previous section, then the spins in $L_1$ are independent
one those in $L_2$ are given. From then on the subsets $L_i$ can be constructed
by similarity.
If periodic boundary conditions are imposed on $L_0$, they inherited
by every successive sublattice.

If one wants to carry out an expansion in
a basis with more polynomials, the requirement that the spins in $L_i$ be
independent once those in $L_{i+1}$ are known places restrictions on the
polynomials one can use; for example, the polynomial $\sum s_{i,j}s_{i+2,j}$
for $i,j$ such that $i+j$ is even is a function of only the spins in $L_1$
but it cannot be used, because the spins at points where
$i+j$ is even while each of $i,j$ is odd would not be independent once the
spins in $L_2$ are known, given the linkage created by this added polynomial.
This places creates a restriction on the accuracy of the
various marginals; see also the discussion in section 6. However, the
point made in the present paper is that significant inaccuracy in the
marginals can be tolerated.

The analogous construction in three dimensions is neither so simple nor
unique. Here is what is done in the present paper:
$L_0=L$ consists of all the nodes on a regular cubic lattice, i.e., the
set of points $(i,j,k)$ where $i,j,k$ are integers between $1$ and $N$ and $N$ is a power of $2$. 

$L_1$ consists of the points$(i,j,k)$ where $i+j$ is even for $k$ odd and odd when $k$ is even; the neighbors of $(i,j,k)$ in $L_1$ are the $8$ points $(i,j \pm 1,k \pm
1),
(i \pm 1, j \pm1,k)$.

$L_2$ consists of the points $(i,j,k)$ where $i$ is even,  $j$ is odd and $k$ is even.
The neighbors of a point $(i,j,k)$ in $L_2$ are the $12$ points $(i \pm 1, j, k \pm 1), (i \pm 1, j \pm 2, k \pm 1)$.

$L_3$ consists of the points $(i,j,k)$ where $i,j,k$ are all odd;
the neighbors of $(i,j,k)$ are $(i \pm 2,j,k), (i,j \pm 2,k), (i,j,k \pm 2)$.
This sublattice is similar to the original lattice with the distance between
sites increased to $2$; the next lattices up can then be obtained by  similarity.

This process of constructing sublattices with ever smaller numbers of spins stops when
one reaches a number of spins small enough to be sampled directly, i.e., by listing all
the states, evaluating their probabilities, and picking a state with a frequency
equal to its probability. One has to decide what the smallest lattice is; the best one could do in this
sequence is a lattice similar to $L_2$ with the points $(i,j,k)$ where $i-1=2\ell$,
$j-1=2\ell$, and $k=2\ell$ for an integer $\ell$ chosen so that there
are  $16$ points in this smallest lattice. Here too each lattice inherits periodic boundary conditions
from the original lattice; on the smallest lattice one notes that due to periodicity
$i+2\ell=(i-2\ell)mod(N)$ so that some of the neighbors of a points $(i,j,k)$ are not
distinct and this must be reflected in the evaluation of the last Hamiltonian,
or else all the linear systems one solves in the projection step are singular.

The polynomial basis used in this paper, except in the diagnostics sections,
consists at every level of polynomials of the form
$\psi_m=\sum s_{i,j,k} s_{i_m,j_m,k_m}$, where $(i,j,k)$ is a point in
the sublattice $L_m$ and $(i_m,j_m,k_m)$ is one of its near neighbors
on that sublattice.

\section{Sampling strategy and weights.}

If the marginals whose computation has just been described were exact,
the algorithm outlined in the introduction would be exact. However,
the marginals are usually only approximate,
and
one has to take into account the errors in them, due both to the use of too few basis functions and to the errors in the numerical
determination of the projection coefficients.
The idea here is to compensate for these errors through appropriate weights.

Suppose one wants to compute the average of a function $h(S)$ of a random variable $S$, whose pdf is $P(x)$, i.e., compute
$E[f(S)]=\int h(x)P(x)dx$.
Suppose one has no way to sample $S$ but one can sample a nearby variable
$S_0$ whose pdf is $P_0(x)$.
One then writes
\begin{eqnarray*}
\int h(x)P(x)dx&=&\int h(x) \frac{P(x)}{P_0(x)} P_0(x) dx \\
&=& E\left[h(S_0)\frac{P(S_0)}{P_0(S_0)}\right]\cong\frac{1}{N_s}\sum h(S_{0i})w_i,
\end{eqnarray*}
where the $S_{0i}$ are successive samples of $S_0$, $i=1,\ldots,N_s$, and
$w_i=P(S_{0_i})/P_0(S_{0_i})$ are sampling weights
(see e.g. \cite{Li1}).
In our case, $P$ is the true probability density $e^{W^{(0)}}/Z$ and $P_0$ is the probability density of the sample $S_{0i}$ produced  by the algorithm
we describe, whose pdf $P_0$  differs from $P$ because the marginals used are only
approximate.

The probability $P_0$ of a sample $S_0=(s_{i,j}, \dots, s_{N,N})$ has to be computed as the sample is produced:
the probability of each state of the spins in $L_n$ is known and therefore
the probability of the starting sample of $L_n$ is known; each time one
samples a spin in $L_j$, $j<n$, one has choices whose probabilities can be
computed.
As a practical matter, one must keep track not of the probabilities themselves but rather of their
logs, or else one is undone by numerical underflow.
Note that in the evaluation of $P/P_0$ the factor $Z^{-1}$ remains unknown, but
as $Z$ is common to all the samples, this does not matter.
(and this
remark can be made into an effective algorithm for evaluating $Z$ and
hence the entropy). In practice I found it convenient to pick a value for $Z$ so that $E[log(P/P_0)]=0.$

In practice, for lattices that are not very small, there is a significant range
of weights, and there is a danger that the averaging will be dominated by
a few large weights, which increase the statistical error.
This issue has been discussed before (see e.g.\cite{Li1}) where it is suggested
that one resort to ``layering";  this is indeed what we shall do, but more cautiously than suggested in previous work. Suppose one caps all weights at some value $W$, i.e., replace the weights $w_i$
by $w_i'=min(w_i,W)$.
The effective number of samples in the evaluation of the variance is $N_W +\sum (w_i/w)$, where $N_W$ is the number of samples with $w_i \ge W$, and the summation is
over the samples with $w_i<W$.
Define the fraction $f$ as the fraction of the samples such that
$w_i>W$ (so that $w_i'=W$);
as $W$ increases the fraction $f$ tends to zero.
The averages computed by the algorithm here depend on $f$ (or $W$) and so does
the statistical error; one has to ascertain that any
result one claims
is  independent of $f$. Typically, as the size of the lattice increases,
the range of weights increases, and therefore the effective number of samples decreases for a given number of samples $N_s$. What one has to do is check that the results converge to
a limit as $f$ decreases while the number of samples is still large
enough for the results to be statistically significant. This may require
an increase in $N_s$ as $N$ increases.

\section{Bootstrapping the sampling.}

So far it has been assumed that one can sample the density $e^{W^{(0)}}/Z$ well enough to compute the coefficient in the Kadanoff expansion (\ref{expansion}) of the marginals.
However, these coefficients are needed to make the
sampling efficient when it would not otherwise be so, and the sampling has to be ``bootstrapped" by iteration so it can be used to
determine its own coefficients.

First, make a guess about the coefficients in (\ref{expansion} ), say, set
$a_{i,j}^{\ell,m}=a_{i,j}^{\ell,m,0}$ for the $\ell-$the coefficient at the point $i,j$ in the sublattice $L_m$, where the numbers $a_{i,j}^{\ell,m,0}$ are some plausible guesses. My experience is that it does not much matter what these guesses are;
I typically picked them to be some moderate constant independent of $i,j,m,\ell$.
Use these coefficients in a sampling procedure to find new coefficients
$a_{i,j}^{\ell,m,1}$, and repeat as necessary.
An iteration of this kind was discussed in \cite{Ch00}, where it was shown that as the number of samples increases the error in the evaluation can be surprisingly small and that the optimal number of polynomials to use for a given overall
accuracy depends on the number of samples.
In the present work I found by numerical experiment that convergence
is faster if, after one evaluates a new coefficient $a_{i,j}^{\ell,m,r+1}$, one sets in
the next round
$a_{i,j}^{\ell,m,r+1}=(a_{i,j}^{\ell,m,r}+a_{i,}^{\ell,m,r+1})/2$.
I found experimentally that there is no advantage in computing these coefficients very
accurately, indeed  a relatively small number of samples is sufficient for each
iteration, and two iterations have been sufficient for all the runs below.

\section{Example 1: The two-dimensional Ising model.}

To check the algorithm and gauge its performance, we begin by applying it to the two-dimensional
Ising model.
I did not write a special program for this case and did not take advantage of the
simplifications which arise when the coupling constants in the Hamiltonians are independent of location and one could replace
four basis functions by the single function consisting of
their sum, and the single expansion coefficient is the same at all points so that the projection can be averaged in space as well as over samples.
Once expansion coefficients have been determined, they can be
used to generate as many samples as one wants.

First, I computed the mean magnetization $E[\mu]$, where $\mu=\frac{1}{N^2}\sum s_{i,j}$, $1\leq i,j \leq N$, as a function of the temperature $T$.
To calculate such means near $T=T_c$ one needs a way to break the symmetry of
the problem; this is usually done by adding a small asymmetric term 
$\eps\sum s_{i,j}$ to the Hamiltonian, for example with $\eps=\eps_0/N$, $\eps_0\sim0.2$.
Adding such a field here works for small $N$, but as $N$ increases,
a value of $\eps_0$ small enough not to disturb the final result may not
suffice to bias the smallest lattice $L_n$ in one direction.
The remedy is to assign positive weights to only to those spin configurations
in $L_n$ (where weights are explicitly known) such that $\sum_{L_n} s_{i,j}\geq0$.

It may be tempting to introduce a symmetry breaker into the initial Hamiltonian
$W^{(0)}$, add terms odd in the $s_{i,j}$ to the series of linkages, and attempt
to compute appropriate symmetry-breaking terms for the smaller lattices
by a construction like the one above. This is not a good idea. The longer series
is expensive to use and the computation of higher Hamiltonians is unstable
to asymmetric perturbations, generating unnecessary errors.

In Table 1 we present numerical results for $T=2.2$ and different values of $N$
compared with values obtained by Metropolis sampling with many sweeps of the lattice.
$N_s=
1000$ samples were used in each of two iterations to calculate the approximate
marginals; once these are found  one can inexpensively generate as many samples of
of the spins as wanted; here 1000 were used.
The Table exhibits the dependence of the computed $E[\mu]$ on the fraction $f$
of weights which have been capped; the statistical error in the estimates of
$E[\mu]$ grows as $f$ decreases, but more slowly than one would expect.
For $N=16, 32, 64$ the results converge as $f\rightarrow 0$ before the statistical error
becomes large,
but not when $N=128$, and one should conclude that this value of $N$ is too large for the present
algorithm with so small a basis in the computation of marginals. Even with
$N=128$ one obtains a reasonable average $(E[\mu]=.80)$ if one is willing to use enough 
samples.
Note also that the weights become large, and the calculations must be performed
in double precision.

\begin{table}[!t]
\begin{center}
\begin{tabular}{ccrlrc}
\multicolumn{6}{c}{\bf Table 1}\\
\\
\multicolumn{6}{c}{Ising magnetization at $T=2.2$}
\\
\\
\hline
size of & no. of samples & $\log W$ & $f$ & $E[\mu]$ & metropolis\\
array      & $N_s$ & & & \\
\hline
$16\times16$&1000&2&.33&.74$\pm$0.01&$.01\pm.01$\\
&&4&.08&.80$\pm$0.01\\
&&6&.00&.80$\pm$0.01\\
\hline
$32\times32$&1000&5&.17&.76$\pm$0.01&$.81\pm.01$\\
&&7&.10&.80$\pm$0.01\\
&&9&.04&.81$\pm$0.015\\
\hline
$64\times64$&1000&15&.134&.74$\pm$.01&$.801\pm.001$\\
&&20&.042&.77$\pm$.01\\
&&25&.009&.80$\pm$.01\\
&&30&.001&.80$\pm$.015\\
\hline
$128\times128$&1000&25&.074&.67$\pm$.01&$.798\pm.001$\\
&&35&.023&.70$\pm$.01\\
&&45&.002&.74$\pm$.02\\
&&50&.001&.75$\pm$.05\\
\multicolumn{6}{c}{no convergence}\\
\hline \end{tabular}
\end{center}
\end{table}

We now turn to the determination of the critical temperature $T_c$.
This can be obtained from the intersection of the graphs of $E[\mu]$ vs.\ $T$
for various values of $N$ (see e.g. \cite{Li1},\cite{Ka3}); it is more instructive here to
apply the construction in \cite{Ch9}, based on the fact that if one expands
the ``renormalized" Hamiltonians in successive linkages, ie., if one finds
the functions $W^{(i)}$ such that the marginals on $L_i$ are $\exp(W^{(i)})/Z$,
using the series (\ref{expansion}), then the coefficients in the series increase when $T<T_c$
and decrease when $T > T_c$.
For this construction to work, one needs enough polynomials for convergence,
i.e., so that the addition of more polynomials leaves the calculation
unchanged.
In the present case, this is achieved with the following 7 polynomials:
$\psi_1,\psi_2,\psi_3,\psi_4=\sum s_{i,j}s_{i\pm1,j\pm1}$ (as above),
$\psi_5=\sum s_{i,j}(s_{i+2,j}+s_{i,j+2}+s_{i-2,j}+s_{i,j-2}),
\psi_6=\sum s_{i,j}\sigma^3_{i,j}/10,
\psi_7=\sum s_{i,j}\sigma^5_{i,j}/100,$
where $\sigma_{i,j}=s_{i+1,j}+s_{i,j+1}+s_{i,j+1}+s_{i,j-1}$.
The use of polynomials with higher powers of the $s_{i,j}$ is essential
(see e.g.\cite{Bin}), and it is the more surprising that the approximate marginals
calculated without them are already able to produce usable samples.
The constant divisors in $\psi_6, \psi_7$ are there to keep all the coefficients within the same
order of magnitude.
No advantage is taken here of the symmetries of the Ising model.
In Table 2 I present the sums of the coefficients in the expansion as a function of the temperature $T$ for the levels $i=2,4,6$ (where the lattices
are mutually similar) in a $N^2$ lattice with $N=16$ ($N$ is chosen small for
reference in the next section). From Table 2 one can readily
deduce that $ 2.25 < T_c < 2.33$; for the value of $T$ between these two
bounds the sum of the coefficients oscillates as $i$ increases. Taking the average of
the two bounds (which are not optimal) yields $T_c\equiv 2.29$ (the exact value
is $T_c=2.269...$). A more careful analysis improves the result and so does a
larger value of $N$. All the coefficients have the same sign, except occasionally when a
coefficient has a very small absolute value.

\begin{table}[!t]
\begin{center}
\begin{tabular}{ccccccccc}
\multicolumn{9}{c}{\bf Table 2}
\\
\\
\multicolumn{9}{c}{Sums of coefficients of Kadanoff expansion}\\
\multicolumn{9}{c}{as a function of $T$ for Ising model}
\\
\\
\hline
$T$&&&$i=2$&&&$i=4$&&$i=6$\\
\hline
2.20&&&1.63&&&2.16&&2.71\\
2.25&&&1.50&&&1.77&&1.88\\
2.26&&&1.49&&&1.73&&1.50\\
2.28&&&1.45&&&1.65&&1.40\\
2.30&&&1.41&&&1.50&&1.30\\
2.32&&&1.36&&&1.40&&1.13\\
2.33&&&1.35&&&1.34&&1.15\\
2.35&&&1.30&&&1.27&&0.95\\ \hline
\end{tabular}
\end{center}
\end{table}
For the sake of completeness I plotted in Figure 1 a histogram of the logarithms
of the weights $w_i$ for the Ising model with $N=32$ and $10^4$ samples;
the zero of $log w$ is chosen as described above. The cost per sample of an optimized version of this program for $N=32,64$ is
competitive with the cost of a cluster algorithm \cite{Sw3} and is significantly lower
than that of a standard Metropolis sampler. It is not claimed that for the Ising model
the chainless sampler is competitive with a cluster algorithm: as $N$ increases
the complexity of the present algorithm grows because one has to add polynomials and/or
put up with a decrease in the number of effective samples, and one also has to do work to examine
the convergence as $f \rightarrow 0$. The present sampler is meant to be useful when MCMC
is slow, as in the next section.

\begin{figure}[tcb]
\begin{center}
  \includegraphics[width=\textwidth]{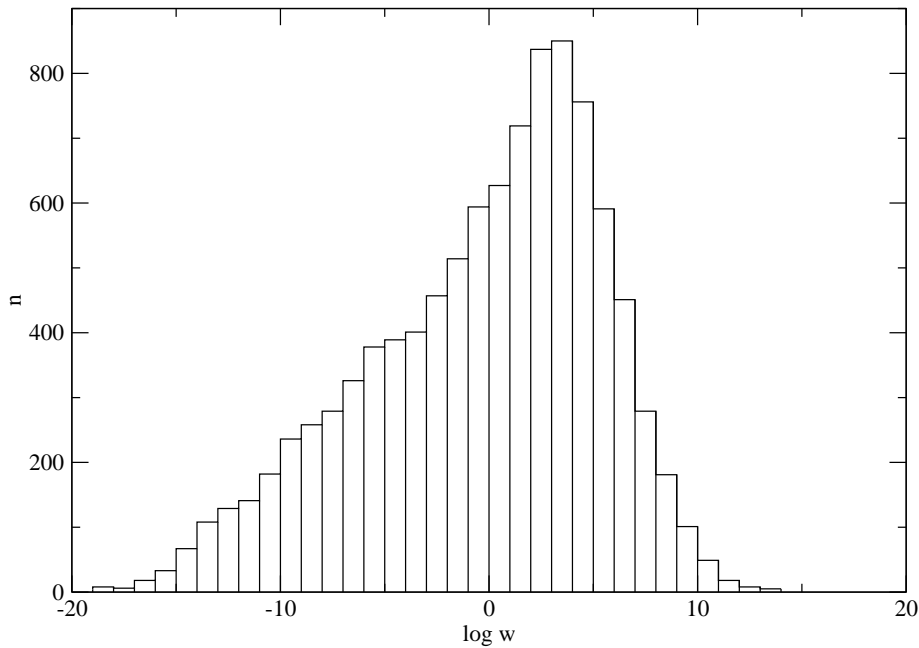}
\end{center}
\caption
{Histogram of weights for the Ising model, $N=32$, $10^4$ samples.}
\label{}
\end{figure}

\section {Example 2: The Edwards-Anderson spin glass in three dimensions}

We now use the chainless construction to calculate some properties of
the EA spin glass \cite{Fi1,Ka3,Me1,Ne3,Ta1}, and in particular, estimate the
critical temperature $T_c$.
The three dimensional Edwards-Anderson spin glass model is defined by equation
(\ref{Hamiltonian}),
where the $J_{i,j,k,\ell}=\beta\xi_{i,j,k,\ell}$, the $\xi_{i,j,k,\ell}$ are independent Gaussian
random variables with mean zero and variance one, and $\beta=1/T$ is the inverse temperature.
Periodic boundary conditions are imposed on the edge of the $N^3$ lattice. 

Let
the symbol $<\cdot>_T$ denotes a thermal average for a given sample of the $J$s, and $[\cdot]_{Av}$
denote an averages over the realizations of the $J$s. Given two independent thermal samples of the
spins $S_{1,2}=\{s^{1,2}_{i,j,k}\}$, we define their overlap to be
$q=N^{-3}\sum s^1_{i,j,k}s^2_{i,j,k}$, where the summation is over all sites in the lattice.
The Binder ratio \cite{Bi3},\cite{Ka3} is $g=0.5(3-[<q^4>_T]_{Av}/[<q^2>_T]^2_{Av})$.
The function $g=g(T)$ is universal, and the graphs of $g$ as a function of $T$ for various
values of the lattice size $N$ should intersect at $T=T_c$. 

The method presented in the present paper is applied to this problem. The only additional
comment needed is that at every point of the lattice one has to invert a matrix generated
by a random process involving integers, and occasionally one of these matrices will be 
singular or nearly so, and will produce unreliable coefficients, particularly for small
samples sizes and low temperatures. As long as there are few such cases, there is no
harm in jettisoning the resulting coefficients and replacing them by zeroes.

In Figure 2 I display the results obtained for this problem. The statistical error is hard to
The numerical parameters are: $2000$ realizations of the $J$s, for each one of them $1000$ samples
for estimating the expansion coefficients and then $5000$ samples for evaluating $q$ and its
moments. I used the bound $log W=30$, which produces a modified
fractions $f=0$ for $N=4$, $f=0.015$ for $N=8$ and $f=.05$ for $N=16$. 
The statistical error was hard to 
gauge; one can readily estimate the standard deviations of the numerical estimates of
$[<q^4>_T]_{Av}$ and $[<q^2>_T]_{Av}$ but these estimates are correlated and one therefore
cannot use their standard deviations to estimate that of $g$. I simply made several runs for some
of these computations and used the scatter of the results to estimate the statistical error.
I concluded that the statistical error is around $1\%$ for $N=4$ and $2-3\%$ for $N=8,16$.

\begin{figure}[tcb]
\begin{center}
  \includegraphics[width=0.95\textwidth]{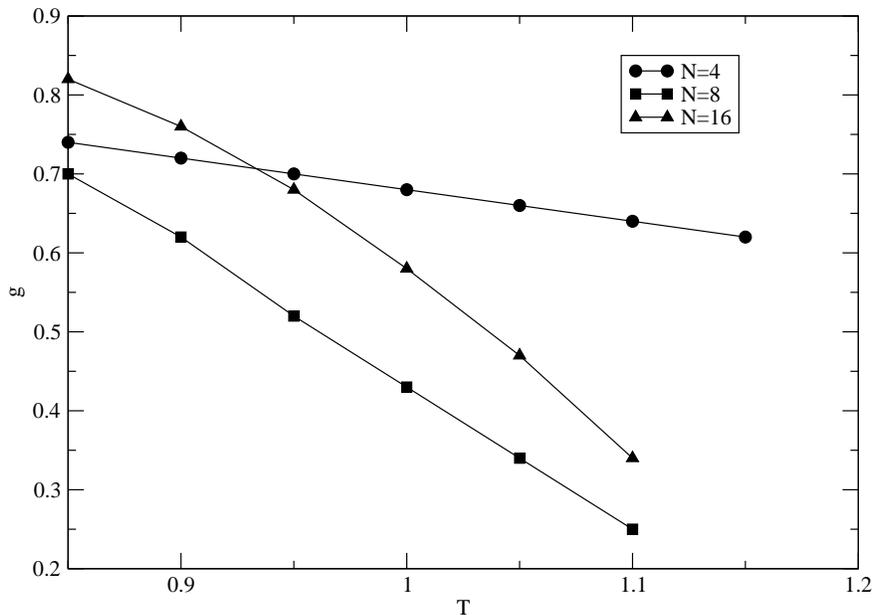}
\end{center}
\caption
{The Binder cumulant $g$ as a function of the temperature $T$ in the three-dimensional AE model.}
\label{}
\end{figure}

These three graphs taken separately approximate the ones in the 
detailed computations of \cite{Ka3}. The graphs for $N=4$ and $N=16$ intersect at $T=.93$,
which is the value of $T_c$ deduced from the Binder cumulant computation in \cite{Ka3}
(and which differs from the value $T_c=.95$ deduced in the same paper from other considerations
and which is
likely to be right from the accumulated previous wisdom, as reported 
in that paper). The graph for $N=8$ is a little off from what one may expect,
but again previous calculations, for example Figure 7 in \cite{Ka3}, also display symptoms
of unexpected waywardness. If one compares Figure 2 with Figure 4 of \cite{Ka3}, one sees
other small discrepancies; for example, the values of $g$ I obtained, in particular for $N=4$, are smaller 
than those in \cite{Ka3} by a small but statistically significant amount; this cannot be
the effect of ``layering" (i.e., the use of a bound $W$) because for $N=4$
layering 
plays no role; it is hard to see how it can be produced the statistical error in either paper
because the sample sizes are certainly large enough. I have no explanation, except the
general orneriness of the EA spin glass, as illustrated by the widely varying results
for various exponents and thresholds summarized in \cite{Ka3}. 
These vagaries do not alter the fact that the algorithm
of the present paper produces worthy results in a very difficult problem.
A more detailed exploration of spin glasses will be published separately. 

It is of interest in the present context to see how the coefficients in the Kadanoff
expansion, used in the previous section to estimate $T_c$ for the Ising model, behave in the
three-dimensional spin glass. Now one needs more polynomials (20; the three-dimensional
analogs of the ones in the preceding section plus others following the same pattern).
The sums of all the coefficients add up to small numbers statistically indistinguishable from
zero, as one may expect, so in Table 3 I display the sums of the absolute values of these coefficients for $N=16$ on self-similar lattices,
which in the Ising case are equally able to exhibit $T_c$ for this value of $N$. No bifurcation between growing
and decreasing sums can be detected near $T_c$, illustrating differences between
the phase transitions in the Ising and spin glass cases. Note the results for $T=0.6$,
a temperature hard to reach by a MCMC process. 

\begin{table}[!t]
\begin{center}
\begin{tabular}{ccccccc}
\multicolumn{7}{c}{\bf Table 3}
\\
\\
\multicolumn{7}{c}{Sums of absolute values of Kadanoff}\\
\multicolumn{7}{c}{coefficients for EA spin glass model}
\\
\\
\hline
$T$&&$i=2$&&$i=5$&&$i=8$\\
\hline
0.6&&3.07&&3.95&&3.73\\
0.9&&2.70&&3.90&&3.42\\
1.0&&2.43&&3.34&&2.76\\
2.00&&6.20&&8.51&&8.87\\ \hline
\end{tabular}
\end{center}
\end{table}

The real question here is the efficiency and speed of the algorithm. What are needed are
timing comparisons between an optimized version of it and optimized version of alternatives,
such as the parallel tempering construction of \cite{Ka3}. This is not available.
The least one can say is that the chainless sampler is highly competitive with others.
Most of the computations in this paper (all but the ones for $g$ at $N=16$) were first
run on a single serial desktop machine.  

\section{Conclusions}

A Monte Carlo sampling technique that relies on a fast marginalization rather than a
Markov chain has been introduced, tested, and applied to a challenging test problem.
The results demonstrate that it is a good alternative, especially for problems where
the free energy has many minima and MCMC algorithms may be slow. Various improvements
to this constructions readily suggest themselves, based on more polynomials, better polynomials,
and
renormalization schemes other than decimation. Related ideas, such as the parallel
marginalization scheme proposed in \cite{We1}, are also worth further investigation
in the context of spin problems. 

The construction of the sequence of lattices above assumed that the original Hamiltonian
involves only near-neighbor interactions; the lifting of this restriction requires
a more elaborate renormalization process and will be pursued elsewhere.

{\bf Acknowledgements} I would like to thank Dr. R. Fattal and Dr. J. Weare for  many illuminating discussions and comments, Profs. E. Mossel and C. Newman for help with the literature, and Dr. M. Lijewski for help in making longer runs. This work was partially supported by the National Science Foundation under
grants DMS-0410110 and DMS-070590, and
by the Director, Office of Science, Computational and Technology Research, U.S. Department of Energy under Contract No. DE-AC02-05CH11231.

\end{document}